\documentclass[12pt,a4paper]{article}
\usepackage{amsmath,amssymb,amsthm,amsfonts}
\usepackage{centernot}
\usepackage{mathtools}
\usepackage{stmaryrd}
\usepackage{color}
\usepackage{enumitem}
\usepackage[english]{babel}
\usepackage{url}
\usepackage{graphicx}
\usepackage{a4}
\newtheorem{theorem}{Theorem}

\newtheorem{prop}[theorem]{Proposition}

\newtheorem{Problem}{Problem}

\usepackage{multirow}

\begin{document}
	\baselineskip 18pt
	
	\title{\bf The lower central series of the unit group of an integral group ring}
	\author{Sugandha Maheshwary 
		\\ {\em \small Indian Institute of Science Education and Research Mohali,}\\
		{\em \small Sector 81, Mohali-140306, Punjab, India.}
		\\{\em \small email: sugandha@iisermohali.ac.in}}

	\date{}
	{\maketitle}
	
	\begin{abstract}The aim of this article is to draw attention towards various natural but unanswered questions related to the lower central series of the unit group of an integral group ring.
	
	\end{abstract}\vspace{.25cm}
	 \vspace{.25cm} 

	\section*{Introduction}

For a group $G$, denote by $\mathcal{V}(\mathbb{Z}G)$, the group of normalized units, i.e., units with augmentation one in the integral group ring $\mathbb{Z}G$. The study of $\mathcal{V}(\mathbb{Z}G)$ and its center attracts a varied set of questions and one naturally seeks the understanding of central series of $\mathcal{V}(\mathbb{Z}G)$. While the upper central series of $\mathcal{V}(\mathbb{Z}G)$ has been well explored, at least for a finite group $G$ (see \cite{MP18} for a detailed survey), apparently, not much is known about its lower central series $\{\gamma_{n}(\mathcal{V})\}_{n\geq 1}$, where $\mathcal{V}:=\mathcal{V}(\mathbb{Z}G)$ and
	$$\gamma_{1}(\mathcal{V})=\mathcal{V},~\gamma_{2}(\mathcal{V})=\mathcal{V}',~\gamma_{i}(\mathcal{V})=[\gamma_{i-1}(\mathcal{V}),\mathcal{V}],~i\geq 2.$$ 

In this article, we try to draw attention towards certain obvious problems \linebreak associated to the study of the lower central series of $\mathcal{V}(\mathbb{Z}G)$, which essentially need to be worked upon.

	\subsection*{The lower central series of $G$ and $\mathcal{V}(\mathbb{Z}G)$}
	
	The structure of $G$ naturally has the implication on $\mathcal{V}(\mathbb{Z}G)$. It is not difficult to see that the lower central series of $G$ and $\mathcal{V}(\mathbb{Z}G)$ coincide, if and only if, $\mathcal{V}(\mathbb{Z}G)=G$. Moreover, for a finite group $G$, this happens precisely if $G$ is an abelian group of exponent $2,\,3,\,4$ or $6$, or  $G=Q_8\times E$, where $E$ denotes an elementary abelian $2$-group and $Q_8$ is the quaternion group of order $8$.

	The  groups for which  $\mathcal{Z}(\mathcal{V}(\mathbb{Z}G))=\mathcal{Z}(G)$ are called {\it cut groups} \cite{BMP17}. A group $G$ is a cut group, if and only if it satisfies the RS property (see \cite{BMP19}). This class of groups possesses various interesting properties and is topic of active research. For instance, for a finite group $G$, being a cut group is  equivalent to saying that the center of $\mathcal{V}(\mathbb{Z}G)$ has rank zero (the rank of $\mathcal{Z}(\mathcal{V}(\mathbb{Z}G))$ has been computed independently by Ferraz \cite{Fer04} and Ritter-Sehgal \cite{RS05}). For more details on finite cut groups, see (\cite{MP18}, Section 3).

	  Analogously, we ask for the classification of groups $G$ for which the derived subgroup of $\mathcal{V}(\mathbb{Z}G)$ coincides with that of $G$.  
	
	\begin{Problem}
		Classify the groups $G$ for which $\mathcal{V}(\mathbb{Z}G)'=G'$.
	\end{Problem} 
	
Clearly, if $\mathcal{V}(\mathbb{Z}G)=G$, then $\mathcal{V}(\mathbb{Z}G)'=G'$.	In fact, the complete classification of finite groups for which $\mathcal{V}(\mathbb{Z}G)'=G'$, is an easy consequence of Hartley-Pickel's result (\cite{HP80}, Theorem 2) seen together with the above stated characterisation of finite groups which have only trivial units in the integral group ring.
	
	\begin{theorem}\label{EqualCommutators}For a finite group $G$, $\mathcal{V}(\mathbb{Z}G)'=G'$, if and only if, $G$ is an abelian group or a Hamiltonion 2-group.\end{theorem}

	If $G$ is an infinite group, the problem remains open. In fact, the classification of groups $G$ such that $ \mathcal{V}(\mathbb{Z}G)=G$ is not known, though it is proved that $ \mathcal{V}(\mathbb{Z}G)=G$, if $G$ is an infinite cyclic group or an ordered group.

	\subsection*{The termination of the lower central series of $\mathcal{V}(\mathbb{Z}G)$}
	
	It is known \cite{AHP93} that for a finite group $G$, the upper central series of $\mathcal{V}(\mathbb{Z}G)$ stabilizes in at  most 2 steps, i.e., the second center equals the third. Analogously, one asks for a group $G$, a bound, if any, on the number of terms in the lower central series of $\mathcal{V}(\mathbb{Z}G)$.
	
	\begin{Problem}
		Given a group $G$, when does lower central series of $\mathcal{V}(\mathbb{Z}G)$ stabilize?
	\end{Problem} 
	
	In general, no bound is known for the number of terms in the lower central series of $\mathcal{V}(\mathbb{Z}G)$. However, if $\mathcal{V}(\mathbb{Z}G)$ is nilpotent, the number of terms in both the upper and the lower central series coincide. For a finite group $G$, it is known that $\mathcal{V}(\mathbb{Z}G)$ is nilpotent, if and only if, $G$ is either abelian or a Hamiltonion 2-group. 
		Further, arbitrary  groups $G$ with  $\mathcal{V}(\mathbb{Z}G)$  nilpotent, are classified.
	
	\begin{theorem}\label{NilpotentInfinite}$(\cite{SZ77a},~\mathrm{Theorem} \,1.1)$
		$\mathcal{V}(\mathbb{Z}G)$ is nilpotent, if and only if, $G$ is \linebreak nilpotent and the torsion subgroup $T$ of $G$ satisfies one of the following conditions:
		
		\begin{description}
			\item[(i)] $T$ is central in $G$.
			\item[(ii)] $T$ is an abelian 2-group and for $x\in G$, $t\in T$, $xtx^{-1}=t^{\pm 1}.$
			\item[(iii)]  $T=E\times Q_{8}$, where $E$ is an elementary abelian 2-group and $Q_{8}$ is the \linebreak quaternion group of order 8. Moreover, $E$ is central in $G$ and conjugation by $x\in G$, induces on $Q_{8}$, one of the four inner automorphisms.
		\end{description}		
	\end{theorem}
	
If $\mathcal{V}(\mathbb{Z}G)$ is not nilpotent, apparently, there is no answer for the stated problem.
	
	\subsection*{The residual nilpotence of $\mathcal{V}(\mathbb{Z}G)$}
	A group $G$ is said to be \emph{residually nilpotent}, if the\emph{ nilpotent residue} defined by $$\gamma_\omega(G):=\cap_{n}\gamma_{n}(G),$$
i.e., the intersection of all members of the lower central series of the group, is trivial.

	As already observed, $\mathcal{V}(\mathbb{Z}G)$ is rarely nilpotent. This is attributed to the presence of non-abelian free groups inside $\mathcal{V}(\mathbb{Z}G)$.  However, a free group is residually nilpotent. Therefore, the possibility of $\mathcal{V}(\mathbb{Z}G)$ being residually nilpotent cannot be ruled out, even when it contains a free subgroup. Consequently, we have the following problem:
	
	\begin{Problem}
		Given a group $G$, when is $\mathcal{V}(\mathbb{Z}G)$ residually nilpotent?
	\end{Problem}
	This problem has also been listed in \cite{Seh93} (Problem 50, p.\,306) and is  neatly answered, in case the the underlying group $G$ is finite.
	
	\begin{theorem} \label{residually_nilpotent}$(\cite{MW82},~\mathrm{Theorem}\, 2.3)$ For a finite group $G$, the group $\mathcal V(\mathbb{Z}G)$ is residually nilpotent, if and only if, $G$ is a nilpotent group which is a $p$-abelian group, i.e., the commutator subgroup $G'$ is a $p$-group, for some prime $p$. 
	\end{theorem}
	
	Further, in \cite{MW82}, substantial results on the residual nilpotence of finitely generated nilpotent and $FC$-groups have been established. These results are mainly driven by the residual nilpotence of the augmentation ideals. Observe that if $\Delta^{(k)}(G)$ denotes the Lie power of $\Delta(G)$, defined as  $$\Delta^{(1)}(G)=\Delta(G),\quad \Delta^{(k+1)}(G)=[\Delta(G),\Delta^{(k)}(G)]\mathbb{Z}G,$$ then,  $\Delta^{(\omega)}(G):=\cap_k \Delta^{(k)}(G)=\{0\}$ implies the residual nilpotence of $\mathcal V(\mathbb{Z}G)$. Moreover, since the residual nilpotence of the augmentation ideal $\Delta(G)$ implies that of its Lie powers, clearly $\Delta^\omega(G)=\{0\}$ implies the residual nilpotence of $\mathcal V(\mathbb{Z}G)$.

	In fact, for a finite group $G$, the following equivalences can be easily observed, in view of $(\cite{Pas79},~ \mathrm{Theorem} \,2.13)$ and Theorem \ref{residually_nilpotent}.
	
	\begin{prop}\label{Equiv_Finite}The following statements are equivalent for a finite group $G$:
		
		\begin{description}
			\item[(i)] $\Delta^{(\omega)}(G)=\{0\}$.
			\item[(ii)] $\mathcal V(\mathbb{Z}G)$ is residually nilpotent i.e., $\gamma_{\omega}(\mathcal{V}(\mathbb{Z}G))=\{1\}.$ 
			\item[(iii)] $G$ is a nilpotent $p$-abelian group.
		\end{description}
	\end{prop}

	Some work on the residual nilpotence of $\mathcal V(\mathbb{Z}G)$ can also be found in \cite{MP20}.

		\subsection*{The quotients of the lower central series of  $\mathcal{V}(\mathbb{Z}G)$}

	In the study of the lower central series of $\mathcal V(\mathbb{Z}G)$, one naturally seeks the understanding of its subsequent quotients. 
	
	\begin{Problem}
		Given a group $G$, describe $\gamma_{i}(\mathcal V(\mathbb{Z}G))/\gamma_{i+1}(\mathcal V(\mathbb{Z}G))$, for $i\geq 1$.
	\end{Problem}

For two non-abelian groups of small orders, we have an answer to this  problem.

	\begin{theorem}$(\cite{SGV97}, \cite{SG00})$ For a group $G$, let $\mathcal{V}:=\mathcal{V}(\mathbb{Z}G)$.
		\begin{description}
			\item[(i)] If $G=S_3$, the symmetric group of order $6$, then the quotient $\mathcal{V}/\gamma_{n}(\mathcal{V})$ is isomorphic to the dihedral group of order $2^n$, for $n\geq 2$. Consequently, for every $n\geq 2$, $\gamma_{n}(\mathcal{V})/\gamma_{n+1}(\mathcal{V})$ is isomorphic to the cyclic group of order $2$.
			\item[(ii)]  If $G=A_4$, the alternating group on $4$ elements, then $\mathcal{V}/\mathcal{V}'$ is isomorphic to the cyclic group of order $3$ and $\gamma_{n}(\mathcal{V})=\mathcal{V}'$, for every $n\geq 2$. 
		\end{description}	\end{theorem}	
	
	A partial study of the above problem for a finite group $G$, namely, the study of the quotient  $\gamma_{1}(\mathcal{V}(\mathbb{Z}G))/\gamma_{2}(\mathcal{V}(\mathbb{Z}G))$, i.e, $\mathcal{V}(\mathbb{Z}G)/\mathcal{V}(\mathbb{Z}G)'$ has been taken up in \cite{BMM20}. \\
	
	Apart from these results, no information on the quotients in the lower central series of $\mathcal{V}(\mathbb{Z}G)$ is readily available, even for a finite group $G$.
	 \subsection*{The lower central series of $\mathcal{V}(\mathbb{Z}G)$: a complete description}

	 	 \begin{Problem}
	 	For a group $G$, give a description of the terms in the lower central series of $\mathcal{V}(\mathbb{Z}G)$.
	 \end{Problem}
 
	 This is a fundamental and the most ideal problem, and of course the most challenging one, that one would seek an answer for, in order to understand the lower central series of $\mathcal{V}(\mathbb{Z}G)$.

 Apparently, the question remains unanswered for any considerable class of groups. However, for some particular groups, namely the dihedral groups of order $6$ and $8$, and the alternating group on $4$ elements, some relevant results can be found in \cite{SGV97},  \cite{SG00} and \cite{SG01}.

\par\vspace{1cm}

\begin{center}
{\bf Acknowledgement} 
\end{center}
\noindent The author's research is supported by DST, India  (INSPIRE/04/2017/000897). This research was also supported in part by the International Centre for Theoretical Sciences (ICTS) during a visit for participating in the program- Group Algebras, Representations and Computation (Code: ICTS/Prog-garc2019/10).

\providecommand{\bysame}{\leavevmode\hbox to3em{\hrulefill}\thinspace}
\providecommand{\MR}{\relax\ifhmode\unskip\space\fi MR }
\providecommand{\MRhref}[2]{%
	\href{http://www.ams.org/mathscinet-getitem?mr=#1}{#2}
}
\providecommand{\href}[2]{#2}

\end{document}